\documentclass[10pt]{amsart}
\usepackage{amsmath}
\usepackage{amssymb}
\usepackage{enumerate}
\usepackage{amsbsy}
\usepackage{amsfonts}
\usepackage{color}

\newtheorem{thm}{Theorem}[section]
\newtheorem{lem}[thm]{Lemma}
\newtheorem{prop}[thm]{Proposition}

\numberwithin{equation}{section}

\theoremstyle{remark}

\newcommand{\les}{\lesssim}
\newcommand{\lam}{{\lambda}}
\newcommand{\gam}{{\gamma}}

\newcommand{\vp}{{\varphi}}
\newcommand{\ve}{{\varepsilon}}
\newcommand{\de}{{\delta}}
\newcommand{\al}{{\alpha}}
\newcommand{\ka}{{\kappa}}

\newcommand{\rtwo}{\mathbb{R}^2}
\newcommand{\rot}{\mathbb{R}^{1+2}}

\newcommand{\psipm}{\psi_\pm}
\newcommand{\psizpm}{\psi_{0,\pm}}
\newcommand{\psip}{\psi_{+}}
\newcommand{\psim}{\psi_{-}}
\newcommand{\pipm}{\Pi^{\pm}(D)}
\newcommand{\lams}{\eta_{\#}}
\newcommand{\kaps}{\kappa_{\#}}

\newcommand{\na}{|\nabla|}

\newcommand{\xsb}{X_{\pm}^{s,b}}
\newcommand{\xsbj}{X_{\pm_j}^{s,b}}

\newcommand{\xob}{X_\pm^{0,b}}
\newcommand{\plam}{P_{\lam}}
\newcommand{\pmu}{P_{\mu}}
\newcommand{\brd}{\Lambda}
\newcommand{\ltwo}{L_{t,x}^2}
\newcommand{\lfou}{L_{t,x}^4}

\begin{document}

\title[LWP of 2D Dirac equation]{Local well-posedness of Dirac equations with nonlinearity derived from honeycomb structure in 2 dimensions}

\author[K. Lee]{Kiyeon Lee}
\address{Department of Mathematics, Jeonbuk National University, Jeonju 54896, Republic of Korea}
\email{leeky@jbnu.ac.kr}

\begin{abstract}
The aim of this paper is to show local well-posedness of 2 dimensional Dirac equations with power type and Hartree type nonlinearity derived from honeycomb structure in $H^s$ for $s>\frac78$ and $s>\frac38$, respectively.  We also provide the smoothness failure of flows of Dirac equations.

\end{abstract}


\maketitle

\section{Introduction}

In this paper we consider following two Cauchy problems for massless honeycomb lattice power type Dirac equations($\ell=1$) and  Hartree type Dirac equations$(\ell=2)$:
\begin{eqnarray}\label{maineq}
\left\{ \begin{array}{c}
(\partial_t + \al\cdot D )  \psi =  -i\kappa \mathcal N_\ell(\psi,\psi)\psi \\
\psi(0)=\psi_0
\end{array} \right.
\end{eqnarray}
where $\psi:\mathbb{R}^{1+2} \to \mathbb{C}^2$ is the spinor field represented by a column vector, $\ka$ is constant, $D = -i\nabla$, and $\al = (P_\#\al^1,P_\#\al^2)$  are the Dirac matrices defined by
\begin{align*}
\alpha^1 &= \left(\begin{array}{ll} 0 & \;\;1 \\ 1 & \;\;0 \end{array} \right), \qquad \alpha^2 = \left(\begin{array}{ll} 0 & -i \\ i & \;\;0 \end{array} \right),
\end{align*}
with $P_\# = \left( \begin{array}{cc} \overline{\lams} & 0  \\ 0 & \lams \end{array} \right)$ for honeycomb lattice constant $\lams\neq 0$ arising from nonlinear Schr\"odinger equations(NLS) with honeycomb lattice potentials(see the Section II in \cite{phymod}). The nonlinearities $N_\ell$ are defined by
$$
\mathcal N_1(\psi_1, \psi_2)=  \left( \begin{array}{cc} b_1 \psi_{11} \overline{\psi_{21}} + 2b_2 \psi_{12}\overline{\psi_{22}}   & 0 \\ 0 &  b_1\psi_{12} \overline{\psi_{22}} + 2b_2 \psi_{11} \overline{\psi_{21}} \end{array} \right), 
\quad \mathcal N_2(\psi_1, \psi_2)=\left(|x|^{-1}*(\psi_1^\dagger \psi_2)\right) 
$$
where $\psi_{j1},\psi_{j2}$ are components of $\psi_j$ and the coefficients $b_1,b_2 >0$ which are an amplitude of Bloch waves. The  symbol $*$ denotes convolution operator in $\mathbb{R}^2$ and the $\psi^\dagger$ is the complex conjugate transpose of $\psi$.

Our main equations with the nonlinearity $\mathcal{N}_\ell$ are derived from two dimensional Schr\"odinger equations with honeycomb lattice potential. Its rigorous derivation appears in \cite{phymod}. The honeycomb lattice structure has  appeared in the fabrication of graphene, a mono-crystalline graphitic film in which electrons behave like massless Dirac fermions(see \cite{castro}). Also, the nonlinear optics which model laser beam propagator in particular types of photonic crystals, have the honeycomb structure(see \cite{bapese, hame}).

The equation \eqref{maineq} for $\ell=1$ has the scaling invariance structure in $\dot{H}^{\frac{1}2}$. That is, for $\psi_1$  the solution  to \eqref{maineq} with $\ell=1$, the function $\psi_{1,\lam}$ defined by $\psi_{1,\lam}(t,x) = \lam^\frac{1}{2} \psi_1 (\lam t,\lam x)$  is also the solution to the equation \eqref{maineq} with $\ell=1$ and satisfies that $\|\psi_{1,\lam}(0,\cdot)\|_{\dot{H}^{\frac12}}= \|\psi_1(0,\cdot)\|_{\dot{H}^{\frac12}} $. By this reason, the \eqref{maineq} for $\ell=1$  is said to be mass-supercritical case. Also, since $\|\psi_{2,\lam}(0,\cdot)\|_{L_x^2}= \|\psi_2(0,\cdot)\|_{L_x^2} $ for $\psi_{2,\lam}(t,x) = \lam\psi_2 (\lam t,\lam x)$ where $\psi_2$ is solution to \eqref{maineq} with $\ell=2$, the equation \eqref{maineq} with $\ell=2$ has the scaling invariance structure in $L_x^2$. The \eqref{maineq} for $\ell=2$ is called to be mass-critical case.




Now we state the main theorem of this paper. For simplicity of representation, we set an index $s(\ell)$ by
\begin{align*}
	s(\ell) = \left\{ \begin{array}{cc}
		\frac12    & \;\mbox{ if } \ell=1, \\
		0          & \;\,\mbox{ if } \ell=2.
	\end{array}  \right.
\end{align*}

\begin{thm}[Local well-posedness for $H^s$ data]\label{mainthm1}
Let $s > s(\ell) + \frac38$ for $\ell=1,2$. Then \eqref{maineq}  is locally well-posed for initial data in $H^s(\rtwo)$.
\end{thm}

Here  a definition of the fractional Sobolev space $H^s(\mathbb R^2)$ is placed in \textbf{Notations} below. In particular, LWP result of Dirac equations which have same nonlinearity $\mathcal{N}_\ell$ has been studied  in \cite{phymod} for $H^{s}(\mathbb{R}^2)$ with $s>1$.

We can prove Theorem 1 for massive cases ($m>0$) in the same way as proof of Theorem 1. Since the Physical model comes from massless Dirac Fermions, we only consider the massless case($m=0$) in this paper.

Lemma \ref{bies} is deduced from Selberg's estimates and  we get a coefficient $\mu^{\frac38-s}$ in the result (3) of Lemma \ref{bies}. Then the condition $s> s(\ell) + \frac38$ is necessary in process of proof of Theorem \ref{mainthm1} and the coefficient $\mu^{\frac38-s}$ makes the gap between scaling critical index $s(\ell)$ and our well-posedness index $s(\ell)+\frac38$.

In this paper, we consider Dirac equations with some nonlinearity.  Related equations to \eqref{maineq} are well known as semi-relativistic equations as follows:
\begin{align}
iu_t + \sqrt{m^2 - \Delta} u &= \lam |u|^2u, \label{cubic}\\
iu_t + \sqrt{m^2 - \Delta} u &= \lam (|x|^{-1}*|u|^2)u. \label{hartree}
\end{align}
The Cauchy problem for semi-relativistic equations with power type nonlinearity \eqref{cubic} has been investigated in \cite{ din, fugeooz}. In \cite{din, fugeooz}, Dinh(\cite{din}) showed local well-posedness(LWP) of \eqref{cubic} with massless case($m=0$) for $H^s(\mathbb{R}^2)$ with $s>\frac34$ and Fujiwara, Georgiev, and Ozawa extended LWP to global well-posedness(GWP) for $H^1(\mathbb{R}^2)$. The Cauchy problem for 3 dimensional Hartree type semi-relativistic equations \eqref{hartree} has been investigated in \cite{len,hele}. First the result of well-posedness was obtained by \cite{len} in $H^s(\mathbb{R}^3)$ for $s \ge \frac12$. In \cite{len}, global well-posedness holds in $H^{\frac12}(\mathbb{R}^3)$ for small data in $L_x^2$. Later this was improved to $s>\frac14$ in  \cite{hele}. Also they(\cite{hele}) showed  ill-posedness result for $H^s(\mathbb{R}^3)$ with $s<\frac14$. For \eqref{hartree} with $d$-dimensions($d\ge2$), Cho and Ozawa(\cite{choz}) have revealed the Global well-posedness result for $H^{s}(\mathbb{R}^d)$ with $s \ge \frac12$.  Further results for semi-relativistic equations, we refer to \cite{chozsashim}.

The difficulty stems from the absence of null-structure of $\mathcal{N}_\ell$. We describe the difference between $\psi^\dagger \beta \psi$ and  $|\psi|^2$ where $\beta = \left(\begin{array}{ll} 1 & \;\,\,\,0 \\ 0 & -1 \end{array} \right)$. The quadratic term $\psi^\dagger \beta \psi$ has a null-structure which represents like $|\psi_1|^2 - |\psi_2|^2$. On the other hand, another term $|\psi|^2 = |\psi_1|^2 + |\psi_2|^2$ does not have the null-structure.  Since this structure induces delicate bilinear estimates, Dirac equations with null-structure lead to better  results than the case without null-structure.  However,  we do not use this structure, because our nonlinearities $\mathcal N_1$ are essentially the same as $|\psi|^2$. For this reason, it is picky to control the nonlinear term $\mathcal N_\ell$. Hence we describe the Lemma \ref{bies} which used crucially in the proof of Theorem \ref{mainthm1}.

Also  we consider the Dirac equation with Coulomb type nonlinearity which has null-structure:
\begin{align}
(i\partial_t + \al \cdot D)\psi &= \lam \left(\psi^\dagger \beta \psi  \right) \psi,\label{powertype}\\
(i\partial_t + \al \cdot D)\psi &= \lam \Big(|x|^{-\gam}*\left(\psi^\dagger \beta \psi \right) \Big)\psi.\label{coulomb}
\end{align}
As known result for the equation \eqref{powertype}, Bejenaru and Herr(\cite{behe}) showed the GWP in $H^{\frac12}(\mathbb{R}^2)$.    And the known results for the equation \eqref{coulomb} are in \cite{choleeoz, lee}.  In \cite{lee}, An author of this paper revealed the LWP in $H^s(\mathbb{R}^2)$ with $s>{\frac{\gam-1}2}$ and $1 \le \gam <2$.  It was studied in  \cite{choleeoz} that global well-posedness and small data scattering holds in $H^s(\mathbb{R}^2)$ for $s>{{\gam-1}}$ and $1 < \gam <2$.   As related to \eqref{coulomb}, there is a Dirac equation with Yukawa potential. One may find many results of the Dirac equation which has Yukawa potential nonlinearity in \cite{cholee, tes2d, tes3d, yang}.

In view of scaling we  expect that  LWP results for \eqref{maineq} is optimal in $H^{s(\ell)}$. For this expectation we introduce the following theorem which denotes the smooth failure of our main equation \eqref{maineq} for $s< s(\ell)$. 
\begin{thm}\label{mainthm2}
	 Let $s < s(\ell)$ and $T>0$. If the flow map $\phi \mapsto u$ in \eqref{maineq} exists as a map from $H^s(\rtwo)$ to $C([-T,T];H^s(\rtwo)),$ it fails to be $C^3$ at the origin.
\end{thm}

If the equation \eqref{maineq} has well-posedness in $[-T,T]$ for some $T>0$, the flows of \eqref{maineq} have the smoothness in $[-T,T]$. Since Theorem \ref{mainthm2} implies that the smoothness of flows of \eqref{maineq} fails, This yields the ill-posedness of \eqref{maineq} for $H^s$ with $s< s(\ell)$.

The smoothness failure of some equations was studied for many authors in \cite{mst, beta, hele, lee}. Molinet, Saut, and Tzvetkov(\cite{mst}), Bejenaru and Tao(\cite{beta}), and Herr and Lenzmann(\cite{hele}) have proved the ill-posedness results similar to Theorem \eqref{mainthm2} for Benjamin-Ono equations, 1-d Schr\"odinger equations and semi-relativistic equations, respectively. For Dirac equation, ill-posedness results have been shown  in \cite{ lee}.

It is still opened the well-posedness of \eqref{maineq} in $H^s(\mathbb{R}^2)$ for $s(\ell) \le s \le s(\ell)+\frac38$. For filling up this gap, we have to obtain better bilinear estimates than Lemma \ref{bies}. For this purpose we should find some structure of nonlinearity of \eqref{maineq} like null-structure. Then we may improve LWP in $H^s$ with sobolev index $s$ below $s(\ell)+\frac38$.

The paper is organized as follows: In Section 2, we discuss  projection operators. In Section 3, we introduce the function spaces and the bilinear estimates the most useful on proof of main theorem. Section 4, we prove Theorem \ref{mainthm1} via the standard contraction method. In section 5, we establish the proof of  Proposition \ref{prop-contract} arising in Section 4.   In the last section, we discuss Theorem \ref{mainthm2} by contradiction argument.

\noindent\textbf{Notations.}\\
\noindent$\bullet$ Space and space-time Fourier transform: $\widehat{f} = \mathcal F_x( f)$ denotes the space variable Fourier transform of $f$ and $\mathcal F_{\xi}^{-1}(g)$ the inverse Fourier transform of $g$ such that
$$
\mathcal{F}_x(f)(\xi) = \int_{\mathbb{R}^2} e^{- ix\cdot \xi} f(x)\,dx,\quad \mathcal F_{\xi}^{-1} (g)(x) = (2\pi)^{-2}\int_{\mathbb{R}^2} e^{ix\cdot \xi} g(\xi)\,d\xi.
$$
$\widetilde{f} = \mathcal F_{t,x}(f)$ denotes the space-time variables Fourier transform of $f$ such that
$$
\mathcal{F}_{t,x}(f)(\tau,\xi) = \int_{\mathbb{R}^{1+2}} e^{-it\cdot\tau- ix\cdot \xi} f(t,x)\,dtdx.
$$

\noindent$\bullet$ Fractional derivatives and Sobolev spaces: $D^s = (-\Delta)^\frac{s}2 = \mathcal F_x^{-1}|\xi|^s\mathcal F_x$, $\Lambda^s = (1- \Delta)^\frac{s}2 = \mathcal F_x^{-1}(1+|\xi|^2)^\frac s2\mathcal F_x$ for $s > 0$. Let us denote $  \dot{H}^s = D^{s}L^2 $ and $H^s = \Lambda^{s}L^2 $ for $s\in \mathbb{R}$.

\noindent$\bullet$ Mixed-normed spaces: For a Banach space $X$ and an interval $I$, $u \in L_I^q X \cap \mathbb{C}$ iff $u(t) \in X$ for a.e.$t \in I$ and $\|u\|_{L_I^qX} := \|\|u(t)\|_X\|_{L_I^q} < \infty$. Especially, we denote  $L_I^qL_x^r = L_t^q(I; L_x^r(\rtwo))$, $L_{I, x}^q = L_I^qL_x^q$, $L_t^qL_x^r = L_{\mathbb R}^qL_x^r$. For vector-valued function $\psi \in L_I^q X \cap \mathbb{C}^2$, we also denote that $\|\psi \|_{L_I^qX} := \| |\psi| \|_{L_I^qX}$.

\noindent$\bullet$ Littlewood-Paley operators: Let us define $\beta_1 \in C^\infty_0(-2,2)$ such that $\beta_1(s) =1$ if $|s|\le 1$ and $\beta_{\lam}(s):= \beta(\frac s\lam) - \beta(\frac{2s}\lam)$ for $\lam >1$. Then we define the frequency projection $\mathcal{F}(\plam f)(\xi) = \beta_{\lam}(\xi)\widehat{f}(\xi)$, $P_{\le \lam} = \sum_{\mu =1}^{\lam}\pmu$ and $P_{\ge \lam} = I- P_{\le \frac{\lam}2}$. Also, for measurable set $S \subset \rtwo, R \subset \mathbb{R}^{1+2}$, we denote that $\mathcal{F}_x(P_{S} f)(\xi) = \chi_{S}(\xi)\widehat{f}(\xi)$ and $\mathcal{F}_{t,x}(P_{R} f)(\tau,\xi) = \chi_{R}(\tau,\xi)\widetilde{f}(\tau,\xi)$.

\noindent$\bullet$ As usual different positive constants depending only on $\al, \ka$ are denoted by the same letter $C$, if not specified. $A \lesssim B$ and $A \gtrsim B$ means that $A \le CB$ and
$A \ge C^{-1}B$, respectively for some $C>0$. $A \sim B$ means that $A \lesssim B$ and $A \gtrsim B$.



\section{Preliminaries}

In this section, for simplicity of the Cauchy problem, we define the projection operators and rewrite the equations \eqref{maineq} to integral equations.

\subsection{Projection operator}
We first define the projections about \eqref{maineq} as following:
\begin{align*}
\Pi^{\pm}(D) := \frac12 \left(I\pm \frac{\al\cdot D}{ |\lams| |\nabla|}  \right).
\end{align*}
Then we get
$$
\al \cdot D = |\lams||\nabla| \left(\Pi^{+}(D) - \Pi^{-}(D) \right).
$$
Using these projection operators, we decompose
$$
\psi = \psip + \psim
$$
where $\psipm :=\pipm\psi$. Also, these projection operators satisfy that
$$
\pipm\pipm=\pipm, \qquad \pipm\Pi^{\mp}(D)=0.
$$
Applying these operator to \eqref{maineq} we see that
\begin{eqnarray}\label{releq}
( \partial_t \pm |\lams||\nabla|) \psipm = -i \kappa \Pi^{\pm}(D)  \mathcal N_\ell(\psi,\psi)\psi
\end{eqnarray}
for $\ell=1,2$ with initial data
$$
\psipm(0) =: \psi_{0,\pm} \in H^s.
$$

To simplify the representation of \eqref{releq},  we set the spinner $\phi_\pm(t,x)= \psipm \left(\frac{t}{|\lams|},x \right)$. Hence $\phi$ satisfies that
$$
(i \partial_t \pm |\nabla|) \phi_\pm = \frac{1}{|\lams|}(i \partial_t \pm |\lams||\nabla|) \psi_\pm =  -\frac{i\kappa}{|\lams|}\Pi^{\pm}(D) \mathcal N_\ell(\phi,\phi)\phi
$$
for  $\ell=1,2$. We still call the spinner to $\psi$. Then we finally get the second main equation
\begin{eqnarray}\label{releq2}
(i \partial_t \pm |\nabla|) \psipm = -i\kaps \Pi^{\pm}(D) \mathcal N_\ell (\psi,\psi)\psi.
\end{eqnarray}
where $\kaps = \frac{\kappa}{|\lams|}$.

By Duhamel's formula, we can represent the equations \eqref{releq2} written as an integral equation
\begin{eqnarray}\label{duhamel}
\psipm(t) &=& S_\pm (t)\psizpm + \kaps\int_{0}^{t}S_\pm(t-t')\pipm \left[\mathcal N_\ell \Big(\psi(t'),\psi(t')\Big)\psi(t')\right] dt'
\end{eqnarray}
for $\ell=1,2$. Here we define the linear propagator $S_\pm( t)$ as following:
\begin{eqnarray}\label{propa}
S_\pm (t)f = e^{\mp it|\nabla|}f.
\end{eqnarray}

\subsection{Fractional Leibniz rule}
The following lemma which is called fractional Leibniz rule is useful in the proof of LWP. 
\begin{lem}[\cite{kapo,kpv2,kgo}]\label{leib}
Let $0 < s < 1,\; 1<p<\infty$. Then
$$
\|D^s(fg)-fD^sg-gD^sf\|_{L^p} \les \|D^{s_1}f\|_{L^{p_1}}\|D^{s_2}g\|_{L^{p_2}}
$$
provided $ s= s_1 + s_2 $ with $0 \le  s_1,s_2 \le 1$ and $\frac1p= \frac{1}{p_1} + \frac{1}{p_2}$.
\end{lem}
The proof of Lemma \ref{leib} is in \cite{kapo,kpv2,kgo}.

\section{Function spaces and Bilinear estimates}

\subsection{Functions spaces: $X^{s,b}$-space}
We first introduce $\xsb$ space which will be useful in local theories. (See e.g. \cite{bourgain, kpv, tao}.)
Let us define the norm for $s,b \in \mathbb{R}$ as follows:
$$
	X_{\pm}^{s,b}(T) := \Big\{ \psi : \Big\|\chi_{[-T,T]} \psi \Big\|_{X_{\pm}^{s,b}} < \infty  \Big\}
$$
with a norm
$$
\|\psi\|_{X_{\pm}^{s,b}}:= \left(\int_{\mathbb{R}^{1+2}} \left|\left<\xi\right>^s \left<\tau \mp |\xi|\right>^b \widetilde{\psi}(\tau,\xi) \right|^2d\tau d\xi \right)^{\frac12}.
$$
In particular,  we denote that $X_{\pm_j}^{s,b}$ is $X_{+}^{s,b}$ for $\pm_j =+$ and  $X_{\pm_j}^{s,b}$ is $X_{-}^{s,b}$ for $\pm_j =-$.
These function spaces satisfy the embedding for $b>\frac12$
$$
X_{\pm}^{s,b}(T) \hookrightarrow C([-T,T]; H^s).
$$

\subsection{Bilinear estimates}
\begin{lem}[Theorem 2.1 of \cite{sel}]\label{selb}
Let $\lam>0$ and $L\ge1$.  Let us define the thickened cones 
\begin{align*}
K_{\lam,L}^{\pm}= \Big\{ (\tau,\xi) : |\xi| \les \lam,\; \tau \mp |\xi| = O(L) \Big\} 
\end{align*}
Then
$$
\|P_{K_{\lam,L}^{\pm} \cap (\mathbb{R} \times B_{\mu})}u\|_{L_{t,x}^4} \les \mu^\frac14 \lam^\frac18 L^\frac38 \|P_{K_{\lam,L}^{\pm} \cap (\mathbb{R} \times B_{\mu})}u\|_{\ltwo}
$$
for $u:\rot \to \mathbb{C}$ and any ball $B_{\mu} \subset \rtwo$ with radius $\mu>0$.
\end{lem}
The following lemma is readily obtained by Lemma \ref{selb}.

\begin{lem}\label{bies}\label{bilinear}   Let $s > \frac38$, $b>\frac12$, and $u:\rot \to \mathbb{C}$. Then the following holds:
\begin{enumerate}
\item[$(1)$] $\|P_{B_{\mu}}u\|_{\lfou} \les \mu^{\frac14}\|P_{B_\mu} u\|_{X_{\pm}^{\frac18,b}}$ for $u \in X_{\pm}^{\frac18,b}$ and any ball $B_\mu$ with radius $\mu>0$,
\item[$(2)$] $\|u\|_{\lfou} \les \|u\|_{X_{\pm}^{\frac38,b}}$ for $u \in X_{\pm}^{\frac38,b}$,
\item[$(3)$] $\|\pmu(u_1\overline {u_2})\|_{\ltwo} \les \mu^{\frac38-s}\|u_1\|_{ X_{\pm_1}^{s,b}}\|u_2\|_{ X_{\pm_2}^{s,b}}$ for $\mu>0$,  $u_j \in X_{\pm_j}^{s,b}$.
\end{enumerate}
\end{lem}

\begin{proof} We first prove (1). Lemma \ref{selb} yields that, for $\lam \ge 1$,
\begin{align*}
&\|P_{B_\mu}\plam u\|_{\lfou} \\
& \les \sum_{L\ge1} \mu^\frac14 \lam^{\frac18}L^{\frac38}\|P_{K_{\lam,L}^{\pm} \cap (\mathbb{R} \times B_\mu)} u\|_{\ltwo} \\
&\les \sum_{L\ge1} \mu^\frac14\lam^{\frac18}L^{\frac38}\|P_{K_{\lam,L}^{\pm} \cap (\mathbb{R} \times B_\mu)} u\|_{\ltwo} \les \sum_{L\ge1} \mu^\frac14\lam^{\frac18}L^{\frac38-b}\|L^b \chi_{K_{\lam,L}^{\pm} \cap (\mathbb{R} \times B_\mu)}\widetilde{u}\|_{L_{\tau,\xi}^2}\\
                             &\les \sum_{L\ge1} \mu^\frac14\lam^{\frac18}L^{\frac38-b}\|\left<\tau \mp |\xi|\right>^b P_{B_\mu} \widetilde{u}\|_{L_{\tau,\xi}^2} \les  \mu^{\frac14}\lam^{\frac18}\|P_{B_\mu} u\|_{\xob}.
\end{align*}
Then we have
$$
\|P_{B_{\mu}}u\|_{\lfou} \les \sum_{\lam \ge 1} \|P_{B_\mu}\plam u\|_{\lfou} \les \sum_{\lam \ge 1} \mu^{\frac14}\lam^\frac18 \|P_{B_\mu} P_{\lam}u\|_{\xob} \les \mu^{\frac14}\|P_{B_\mu}u\|_{X^{\frac18,b}_\pm}. 
$$

For (2), by (1) we obtain
\begin{align*}
\|u\|_{L_{t,x}^4} \les \left(  \sum_{\mu\ge1}  \|P_{\mu}u\|_{L_{t,x}^4}^2\right)^{\frac12} \les \left(  \sum_{\mu\ge1}   \mu^{\frac14}\|P_{\mu}u\|_{X_{\pm}^{\frac18,b}}^2\right)^{\frac12} \les \|u\|_{X_{\pm}^{\frac38,b}}.
\end{align*}

Let us now prove (3).  Using frequency localization and (2), we see that
\begin{align*}
	\left\| \pmu (u_1 \overline{u_2}) \right\|_{L_{t,x}^2} &\les \sum_{\begin{subarray}{l} \lam_1,\lam_2 \ge 1 \\ \mu \les \lam_1 \sim \lam_2\end{subarray}} \left\| \pmu (u_{1,\lam_1} \overline{u_{2,\lam_2}}) \right\|_{L_{t,x}^2} + \sum_{\begin{subarray}{l} \;\;\;\lam_1,\lam_2 \ge 1 \\ \lam_{\min} \les \lam_{\max} \sim \mu \end{subarray}} \left\| \pmu (u_{1,\lam_1} \overline{u_{2,\lam_2}}) \right\|_{L_{t,x}^2}\\
	& \les \sum_{\begin{subarray}{l} \lam_1,\lam_2 \ge 1 \\ \mu \les \lam_1 \sim \lam_2\end{subarray}} \|u_{1,\lam_1}\|_{L_{t,x}^4} \|u_{2,\lam_2} \|_{L_{t,x}^4} + \sum_{\begin{subarray}{l} \;\;\;\lam_1,\lam_2 \ge 1 \\ \lam_{\min} \les \lam_{\max} \sim \mu \end{subarray}} \|u_{1,\lam_1}\|_{L_{t,x}^4} \|u_{2,\lam_2} \|_{L_{t,x}^4}\\
& \les \sum_{\begin{subarray}{l} \lam_1,\lam_2 \ge 1 \\ \mu \les \lam_1 \sim \lam_2\end{subarray}} \lam_1^{\frac38-s}\lam_2^{\frac38-s}\|u_{1,\lam_1}\|_{X_{\pm_1}^{s,b}} \|u_{2,\lam_2} \|_{X_{\pm_2}^{s,b}} + \sum_{\begin{subarray}{l} \;\;\;\lam_1,\lam_2 \ge 1 \\ \lam_{\min} \les \lam_{\max} \sim \mu \end{subarray}} \lam_{\max}^{\frac38-s}\|u_{1,\lam_1}\|_{X_{\pm_1}^{s,b}} \|u_{2,\lam_2} \|_{X_{\pm_2}^{s,b}}\\
&\les \mu^{\frac38 -s} \|u_1\|_{X_{\pm_1}^{s,b}}\|u_2\|_{X_{\pm_2}^{s,b}}.
\end{align*}
Here we used $\lam_{\max}= \max(\lam_1, \lam_2)$, $\lam_{\min}= \min(\lam_1, \lam_2)$.
\end{proof}

\section{Local Well-posedness: Proof of Theorem \ref{mainthm1}}

Let us define a complete Banach metric space $(\mathcal M^{s,b}(T,\de),d)$ as follows:
\begin{align*}
\mathcal{M}^{s,b}(T,\de) &= \left\{ \psi \in C\left( [-T,T]: L_x^2\right) \cap X_{\pm}^{s,b}(T) : \|\psi\|_{\mathcal{M}^{s,b}} : = \|\psi_+\|_{X_{+}^{s,b}} + \|\psi_-\|_{X_{+}^{s,b}} <\de   \right\},\\
d(\psi, \phi) &= \|\psi - \phi\|_{\mathcal{M}^{s,b}}.
\end{align*}
We now consider a map $\mathcal{D}$ on $\mathcal{M}^{s,b}(T,\de)$ by
\begin{align*}
\mathcal{D}(\psi) = \sum_{\pm_0 \in \{\pm \}} S_{\pm_0}( t) \psi_{0,\pm_0} +  \sum_{\begin{subarray}{l}
	\pm_j \in \{ \pm\}\\j=0,1,2,3  \end{subarray}} \kappa_{\#} \int_0^t S_{\pm_0}(t-t') \Pi^{\pm_0}(D) \left[\mathcal N_\ell(\psi_{\pm_1},\psi_{\pm_2})\psi_{\pm_3} \right] dt'.
\end{align*}
where  $\sum_{\pm_0 \in \{\pm\}}F_{\pm_0} = F_{+} + F_-$. Then we first show the map $\mathcal{D}$ is self-mapping on $\mathcal{M}^{s,b}(T, \de)$. By Lemma 2.1 of \cite{gtv} we see that 
$$
\Big\|\chi_{[-T,T]}S_{\pm}( t)\psi_{0,\pm} \Big\|_{\xsb} \les T^{\frac12 -b}\|\psi_0\|_{H^s}
$$
and
$$
\left\|\chi_{[-T,T]}\int_0^t S_{\pm}(t-t')f(t')dt' \right\|_{\xsb} \les T^{1-b+b'}\|f\|_{X^{s,b'}}
$$
for $-\frac12 < b' < 0 <\frac12 < b \le b'+ 1$.


\begin{prop}\label{prop-contract}
	Let $s>s(\ell) + \frac38$ for $\ell=1,2$. Then there exists $-\frac12 < b' < -\frac14 <\frac12 < b \le b'+ 1$ and $\ve >0$, such that
$$
\Big\| \mathcal N_\ell(\psi_1, \psi_2) \psi_3\Big\|_{X_\pm^{s,b'}}\le T^{\ve}\prod_{j=1}^{3} \|\psi_j\|_{\xsbj}
$$
for all $\psi:\rot \to \mathbb{C}^2$ and $\psi_j \in \xsbj$ with ${\it supp}(\psi_j) \subset \{(t,x) : |t| \le T\}.$
\end{prop}
Proposition \ref{prop-contract} will be proved in the next section. We now assume the validity of Proposition \ref{prop-contract}. Then we estimate
\begin{align*}
\| \mathcal{D} (\psi)\|_{\mathcal{M}^{s,b}} \le C_1\|\psi_0\|_{H^s} + C_2 T^\ve \sum_{\pm} \|\psi_{\pm}\|_{X_{\pm}^{s,b}}^3 \le C_1 \|\psi_0\|_{H^s}  + C_2 T^\ve \de^3. 
\end{align*}
Set $ C_1\|\psi_0\|_{H^s} < \frac \de2$ and choose the time $T$ that satisfies $C_2 T^{\ve} \de^3 < \frac \de 2$. Hence we see that $\| \mathcal{D} (\psi)\|_{\mathcal{M}^{s,b}} < \de$. Therefore $\mathcal{D}$ is self-mapping on $\mathcal{M}^{s,b}(T,\de)$. We now describe the fact that $\mathcal{D}$ is contraction mapping on $\mathcal{M}^{s,b}(T,\de)$: 
\begin{align*}
d(\mathcal{D}(\psi), \mathcal{D}(\phi)) &= \|\mathcal{D}(\psi) - \mathcal{D}(\phi)\|_{\mathcal{M}^{s,b}} \le C \left( \|\psi\|_{\mathcal{M}^{s,b}}^2  + \|\phi\|_{\mathcal{M}^{s,b}}^2 \right) \|\psi - \phi \|_{\mathcal{M}^{s,b}} \\
&\le 2C\de^2 \|\psi - \phi \|_{\mathcal{M}^{s,b}} < \frac12 d(\mathcal{D}(\psi), \mathcal{D}(\phi))
\end{align*}
for $\de$ satisfying that $4C\de^2 < \frac12$. 

Therefore this completes the proof of the local existence and uniqueness of a solution to \eqref{maineq}.

\section{Proof of Proposition \ref{prop-contract}}

\subsection{Proof of Proposition \ref{prop-contract}}

	By duality, it  suffices to prove that
	$$
	I_\ell := \left|\small\iint  \mathcal N_\ell(\psi_1, \psi_2) \psi_3  \brd^s \psi_4^\dagger dtdx\right| \les T^{\ve} \prod_{j=1}^{3} \|\psi_j\|_{\xsbj} \|\psi_4\|_{X_{\pm_4}^{0,-b'}}
	$$
for $\psi_4 \in X_{\pm_4}^{0,-b'}$ and $\ell=1,2$. We set $\psi_j =  \left(\begin{array}{c}\psi_{j1}\\\psi_{j2}  \end{array} \right)$ for $j=1,2,3,4$. Then we have
	\begin{align*}
I_1 &=\left|\small\iint   \left( \begin{array}{c} \left(b_1 \psi_{11} \overline{\psi_{21}} + 2b_2 \psi_{12}\overline{\psi_{22}} \right)\psi_{31} \\ - \left( b_1\psi_{12} \overline{\psi_{22}} + 2b_2 \psi_{11} \overline{\psi_{21}}\right)\psi_{32}  \end{array}  \right)  \Big( \begin{array}{cc} \brd^s \overline{\psi_{41}} &\;\;\brd^s \overline{\psi_{42}} \end{array}\Big) dtdx\right|\\
&= C \sum_{j,k,l\in\{1,2\}} \left|\small\iint \psi_{1j}  \overline{\psi_{2j}} \psi_{3k}\brd^s \overline{\psi_{4l}} dtdx\right|
\end{align*}
and
	\begin{align*}
	I_2   = \left|\small\iint |\nabla|^{-1}(\psi_1^\dagger \psi_2)\psi_3\brd^s \psi_4^\dagger dtdx\right| = C \sum_{j,k \in \{1,2\}} \left|\small\iint |\nabla|^{-1}  \left(\overline{\psi_{1j}}  \psi_{2j} \right) \psi_{3k}\brd^s \overline{\psi_{4k}} dtdx\right|.
\end{align*}
To compute  the terms above, we introduce  $\mathbb{C}$-valued version estimates below which will be proved Section \ref{prooflem}.

\begin{lem}\label{scal-cub}
       The following two estimates hold:
	\begin{enumerate}
		\item[$(i)$] Let $s>s(1)+ \frac38$. Then there exists $-\frac12 < b' < -\frac14  <\frac12 < b \le b'+ 1$ and $\ve >0$, such that
			\begin{align}\label{nonlinear-cubic}
			\left|\small\iint (u_1  \overline{ u_2} )u_3\brd^s  \overline{ u_4} dtdx\right| \les T^{\ve} \prod_{j=1}^{3} \|u_j\|_{X_{\pm_j}^{s,b}} \|u_4\|_{X_{\pm_4}^{0,-b'}}
		\end{align}
	for all $u_j:\rot \to \mathbb{C}$ and $u_j \in \xsbj$ with ${\it supp}(u_j) \subset \{(t,x) : |t| \le T\}.$
	  \item[$(ii)$] Let $s>s(2)+ \frac38$. Then there exists $-\frac12 < b' < -\frac14  <\frac12 < b \le b'+ 1$ and $\ve >0$, such that
	  	\begin{align}\label{nonlinear-hartree}
	  	\left|\small\iint \left[|x|^{-1}*( \overline{u_1} u_2) \right] u_3\brd^s  \overline{ u_4} dtdx\right|    \les T^{\ve}\prod_{j=1}^{3} \|u_j\|_{X_{\pm_j}^{s,b}} \|u_4\|_{X_{\pm_4}^{0,-b'}}
	  \end{align}
	\end{enumerate}

	for all $u_j:\rot \to \mathbb{C}$ and $u_j \in \xsbj$ with ${\it supp}(u_\ell) \subset \{(t,x) : |t| \le T\}.$
\end{lem}

	By Lemma \ref{scal-cub}, we get
	$$
	I_\ell \les T^{\ve}\sum_{j,k=1,2}\|\psi_{1j}\|_{X_{\pm_1}^{s,b}}\|\psi_{2j}\|_{X_{\pm_2}^{s,b}}\|\psi_{3k}\|_{X_{\pm_3}^{s,b}}\|\psi_{4k}\|_{X_{\pm_4}^{0,-b'}} \les T^{\ve} \prod_{j=1}^{3} \|\psi_j\|_{\xsbj} \|\psi_4\|_{X_{\pm_4}^{0,-b'}}.
	$$
	for $\ell=1,2$. It completes the proof of Proposition \ref{prop-contract}.

\subsection{Proof of Lemma \ref{scal-cub}}\label{prooflem}

\begin{proof}[Proof of $(i)$ of Lemma \ref{scal-cub}]
We first set $\frac78 < s \le 1 $. By H\"older inequality and Lemma \ref{leib}, we can split the left-hand side of \eqref{nonlinear-cubic} as follows:
\begin{align*}
\left|\small\iint (u_1  \overline{ u_2} )u_3\brd^s  \overline{ u_4} dtdx\right|  &\le \left| \iint \brd^s (  \overline{ u_1}  u_2  u_3) \overline{u_4} dtdx - \iint \brd^s(\overline{ u_1} u_2 )u_3\overline{ u_4 } dtdx -\iint \overline{ u_1} u_2 \left(\brd^su_3 \right)\overline{ u_4 }dtdx \right|\\
 &\qquad\qquad +\left|\iint \brd^s (\overline{ u_1}  u_2 ) u_3 \overline{ u_4} dtdx\right| +\left|\iint   \overline {u_1} u_2  \left(\brd^su_3 \right) \overline{ u_4} dtdx\right|  \\ 
 & =: J_1^1 + J_1^2 + J_1^3.
 \end{align*}

We first treat the $J_1^1$. By Lemma \ref{leib}, we estimate
\begin{align*}
J_1^1 & \les  \left\|\brd^s (  \overline{ u_1}  u_2  u_3) -  \brd^s(\overline{ u_1} u_2 )u_3 - \overline{ u_1} u_2 \left(\brd^su_3 \right) \right\|_{L_t^\frac43L_x^2}\|u_4\|_{L_t^4L_x^2}\\
&\les   \| \brd^s (u_1 \overline{ u_2})\|_{L_{t,x}^2}\|u_3\|_{L_t^4L_{x}^\infty}  \|u_4\|_{L_t^4L_x^2}
\end{align*}

Let us consider $\| \brd^s (u_1 \overline{ u_2})\|_{L_{t,x}^2}$. Like above estimates, Lemma \ref{leib}  yields that 
\begin{align*}
	\| \brd^s (u_1 \overline{ u_2})\|_{L_{t,x}^2} & \les \Big\|  \brd^s(u_1 \overline{ u_2}) - (\brd^su_1) \overline{ u_2} -u_1 (\brd^s\overline{ u_2})\Big\|_{L_{t,x}^2} + \Big\|  (\brd^su_1) \overline{ u_2}  \Big\|_{L_{t,x}^2} + \Big\|  u_1 (\brd^s\overline{ u_2})\Big\|_{L_{t,x}^2}\\
	& \les \Big\|  \brd^su_1\Big\|_{L_t^4L_x^2} \|u_2 \|_{L_{t}^4L_x^\infty} +  \|u_1 \|_{L_{t}^4L_x^\infty}\Big\|  \brd^su_2\Big\|_{L_t^4L_x^2}.
\end{align*}
By Sobolev embedding and Lemma \ref{bilinear} we get
\begin{align}\label{uj-esti}
\Big\|  u_j \Big\|_{L_t^4L_x^\infty} \les \Big\|  \brd^{s-\frac38} u_j\Big\|_{L_{t,x}^4} \les \| u_j\|_{X_{\pm_j}^{s,b}}
\end{align}
for $s>\frac78,\; b>\frac12$, and $j=1,2$. By embedding $X^{0,\frac14} \hookrightarrow L_t^4L_x^2$, the estimate \eqref{uj-esti} leads us that
\begin{align}\label{u12-esti}
\| \brd^s (u_1 \overline{ u_2})\|_{L_{t,x}^2} \les \|u_1\|_{X_{\pm_1}^{s,b}}\|u_2\|_{X_{\pm_2}^{s,b}}.
\end{align}
In particular, by \eqref{uj-esti}, we have
\begin{align}\label{u3-esti}
	\|u_3\|_{L_{t}^4L_x^\infty}\les \|u_3\|_{X_{\pm_3}^{s,b}}.
\end{align}
Using the estimates  \eqref{u12-esti}, \eqref{u3-esti}, and embedding $X^{0,\frac14} \hookrightarrow L_t^4L_x^2$ , we see taht
\begin{align*}
	J_1^1 &\les \|u_1\|_{X_{\pm_1}^{s,b}}\|u_2\|_{X_{\pm_2}^{s,b}}\|u_3\|_{X_{\pm_3}^{s,\frac14}}\|u_4\|_{X_{\pm_4}^{0,\frac14}}\\
	&\les T^{\de} \prod_{j=1}^3{\|u_{j}\|_{X_{\pm_j}^{s,b}}}\|u_4\|_{X_{\pm_4}^{0,-b'}}.
\end{align*}



On the other hand, for $J_1^2,J_1^3$, we obtain
\begin{align}
J_1^2 &\les  \| \brd^s (u_1 \overline{ u_2})\|_{\ltwo}\|u_3\|_{L_t^4L_x^\infty}  \|u_4\|_{L_t^4L_x^2},\label{j12}\\
J_1^3 &\les  \|u_1\|_{L_t^4L_x^\infty}\|u_2\|_{L_t^4L_x^\infty} \|\brd^s u_3\|_{L_t^4L_x^2} \|u_4\|_{L_t^4 L_x^2}. \label{j13}
\end{align}
The estimates for \eqref{j12} are obtained in a similar way to estimates of $J_1^1$.  Hence we consider the  \eqref{j13}. Using (2) of Lemma \ref{bies} and Sobolev embedding,  we estimate
\begin{align}\label{second}
	\begin{aligned}
		\|u_j\|_{L_t^4L_x^\infty}   &\les \|\brd^{s-\frac38}u_j\|_{\lfou} \les \|u_j\|_{\xsbj}, \;\; \mbox{for}\;\; j=1,2,\\
		\|\brd^s u_3\|_{L_t^4L_x^2} &\les \|u_3\|_{X^{s,\frac14}} \les \|u_3\|_{X_{\pm_3}^{s,b}}.
	\end{aligned}
\end{align}
Then the estimate \eqref{second} yields that
$$
J_1^3 \les  T^{\de} \prod_{j=1}^3\|u_j\|_{\xsbj}\|u_4\|_{X_{\pm_4}^{0,-b'}}.
$$
Therefore this completes the proof of  \eqref{nonlinear-cubic}.
\end{proof}

\begin{proof}[Proof of $(ii)$ of Lemma \ref{scal-cub}]
	 The LHS of \eqref{nonlinear-hartree} is bounded by 
	\begin{align*}
			&\left|\small\iint \left[|x|^{-1}*( \overline{u_1} u_2) \right] u_3\brd^s  \overline{ u_4} dtdx\right| \\
			 &\le \left| \iint \brd^s (\na^{-1}(  \overline{ u_1}  u_2 ) u_3) \overline{u_4} dtdx - \iint \brd^s\na^{-1}(\overline{ u_1} u_2 )u_3\overline{ u_4 } dtdx -\iint  \na^{-1}(\overline{ u_1} u_2 ) (\brd^su_3)\overline{ u_4 }dtdx \right|\\
		&\qquad\qquad +\left|\iint \brd^s \na^{-1}(\overline{ u_1}  u_2 ) u_3 \overline{ u_4} dtdx\right| +\left|\iint  \na^{-1}( \overline {u_1} u_2 ) (\brd^su_3) \overline{ u_4} dtdx\right|  \\
		&=: J_2^1 + J_2^2 + J_2^3.
	\end{align*}

We first consider the $J_2^1$. Lemma \ref{leib} yields that
\begin{align*}
	J_2^1 & \les  \left\|\brd^s \left[\na^{-1}(  \overline{ u_1}  u_2 ) u_3 \right] -  \brd^s\na^{-1}(\overline{ u_1} u_2 )u_3 - \na^{-1}(\overline{ u_1} u_2) \left(\brd^su_3 \right) \right\|_{L_t^\frac43L_x^2}\|u_4\|_{L_t^4L_x^2}\\
	&\les   \| \brd^s \na^{-1} (u_1 \overline{ u_2})\|_{L_t^2 L_x^4}\|u_3\|_{L_{t,x}^4}  \|u_4\|_{L_t^4L_x^2}\\
	&\les \sum_{\mu} \| P_{\mu} \brd^s\na^{-1}(u_1 \overline{ u_2})\|_{L_t^2 L_x^4}\|u_3\|_{L_{t,x}^4}  \|u_4\|_{L_t^4L_x^2}.
\end{align*}
Using Hardy-Littlewood-Sobolev and Young's convolution inequality we see that		
\begin{align*}
	\| P_{\le1} \brd^s\na^{-1}(u_1 \overline{ u_2})\|_{L_t^2 L_x^4} & \les \|P_{\le1}(u_1 \overline{ u_2})\|_{L_t^2L_x^{\frac43}} = \|\check{\beta_1} * (u_1 \overline{ u_2})\|_{L_t^2L_x^{\frac43}} \les \|\check{\beta_1}\|_{L_x^\frac43} \|u_1\|_{L_t^4L_x^2}\|u_2\|_{L_t^4L_x^2}\\
	&\les \|u_1\|_{X_{\pm_1}^{0,\frac14}} \|u_2\|_{X_{\pm_2}^{0,\frac14}},
\end{align*}
In third inequality, we used the $\|\check{\beta_1}\|_{L_x^p}\les 1$ for $p>1$. And, by Lemma \ref{bilinear}, we estimate
\begin{align*}
\sum_{\mu\ge2}\| P_{\mu} \brd^s\na^{-1}(u_1 \overline{ u_2})\|_{L_t^2 L_x^4}  & \les \sum_{\mu \ge 2}\mu^{s-1}\|\pmu(\overline{ u_1} u_2 )\|_{L_t^2L_x^4} \les \sum_{\mu \ge 2}\mu^{s-\frac12}\|\pmu(\overline{ u_1} u_2 )\|_{L_{t,x}^2}\\
& \les \sum_{\mu \ge 2}\mu^{-\frac18}\|u_1\|_{X_{\pm_1}^{s,b}}\| u_2\|_{X_{\pm_2}^{s,b}} \les \|u_1\|_{X_{\pm_1}^{s,b}}\|u_2\|_{X_{\pm_2}^{s,b}}.	
\end{align*}
Also, by second estimate of \eqref{second}, we get
$$
J_2^1 \les   \prod_{j=1}^3\|u_j\|_{\xsbj}\|u_4\|_{L_t^4 L_x^2} \les T^\de \prod_{j=1}^3\|u_j\|_{\xsbj}\|u_4\|_{X_{\pm_4}^{0,-b'}}.
$$

Estimates for $J_2^2$ are obtained in almost the same way as estimates for $J_2^1$. Hence it is left to deal with $J_2^3$. By Hardy-Littlewood-Sobolev and Young's convolution inequality, we have
	\begin{align*}
	\left|\iint \na^{-1} P_{\le2}\left(\overline{ u_1} u_2 \right) \left(\brd^{s}u_3 \right)\overline{ u_4} dtdx\right| &\les \|\na^{-1}P_{\le2}( \overline{ u_1} u_2 )\|_{L_t^2L_x^4} \left\|P_{\le 2}\left[(\brd^{s}u_3)\overline{ u_4} \right]\right\|_{L_t^2L_x^{\frac43}}\\
	&\les \|P_{\le 2}( \overline{ u_1} u_2)\|_{L_t^2L_x^\frac43}\|(\brd^{s}u_3)\overline{ u_4} \|_{L_t^2L_x^1}\\
	&\les \|u_1\|_{L_t^4L_x^2}\|u_2\|_{L_t^4L_x^2}\|\brd^{s}u_3\|_{L_t^4L_x^2}\|u_4\|_{L_t^4L_x^2}\\
	&\les \|u_1\|_{X_{\pm_1}^{0,\frac14}}\|u_2\|_{X_{\pm_2}^{0,\frac14}}\|u_3\|_{X_{\pm_3}^{s,\frac14}}\|u_4\|_{X_{\pm_4}^{0,\frac14}}.
\end{align*}
Since there is no contribution of $P_{\le1}(\overline{ u_1} u_2)$, we assume that $P_{\le1}(\overline{ u_1} u_2)=0$. Let us consider the high-frequency part of $J_2^3$. By Lemma \ref{bies} and Bernstein's inequality we estimate
	\begin{align*}
		J_2^3 &\les \sum_{\mu \ge 2}\left|\iint  \na^{-1}P_{\mu}( \overline {u_1} u_2 ) (\brd^su_3) \overline{ u_4} dtdx\right|\\
		&\les \sum_{\mu \ge 2}\mu^{-1}\Big\|\pmu( \overline{ u_{\lam_1}} u_{\lam_2} )\Big\|_{L_t^2L_x^\infty} \Big\|  (\brd^su_3) \overline{ u_4}  \Big\|_{L_t^2 L_x^1}\\
	&\les \sum_{\mu \ge 2}\Big\|\pmu( \overline{ u_{\lam_1}} u_{\lam_2} )\Big\|_{L_{t,x}^2} \Big\|  \brd^su_3 \Big\|_{L_t^4 L_x^2}  \Big\|   u_4  \Big\|_{L_t^4 L_x^2}\\
		&\les \sum_{\mu \ge 2} \mu^{\frac38 -s}\| u_{1}\|_{X_{\pm_1}^{x,b}} \|u_{\lam_2} \|_{X_{\pm_2}^{s,b}} \Big\|  \brd^su_3 \Big\|_{L_t^4 L_x^2}  \Big\|   u_4  \Big\|_{L_t^4 L_x^2}\\
		& \les T^{\delta}\prod_{j=1}\|u_j\|_{\xsbj}\|u_4\|_{X_{\pm_4}^{0,-b'}}.
\end{align*}
Here we used the assumption $s>\frac38$ and $b' < -\frac14$. Therefore this completes the proof of the \eqref{nonlinear-hartree}.
\end{proof}

\section{The proof of Theorem \ref{mainthm2}}\label{proof-ill}
This section aims to show Theorem \ref{mainthm2}. It adopts the argument \cite{mst,hele, lee} to prove of smoothness failure of flows of \eqref{maineq} with cubic and Hartree type nonlinearity. As in proof of \cite{mst, hele}, if the flow map $\psi \mapsto u$ is $C^3$ at the origin from $H^s$ to $C( \left[ 0,T\right);H^s)$, we have \eqref{illpose}. In \cite{mst, hele}, they showed smoothness failure of flows of Benjamin-Ono, semi-relativistic equations, respectively. For the results about Dirac equation, we refer to \cite{ lee}.
  Let us consider the system of equation($\ell=1,2$):
\begin{align}\label{maineq-ill}
\left\{
\begin{array}{l}
(\partial_t + \al\cdot D) \psi = -i\ka \mathcal N_\ell(\psi,\psi)\psi,\\
\psi(0)=\delta\psi_0\in H^s(\mathbb R^2).
\end{array}\right.\end{align}

If the flow is $C^3$ at the origin in $H^s$, then it follows that
\begin{eqnarray}\label{dthree}
\partial_\delta^3\psi(0,t,\cdot)  = 6C\sum_{\pm_j,j=1,2,3,4}\int_0^t S_{\pm_1}(t-t')\Pi^{\pm_1}(D) \left[ \mathcal N_\ell\Big( S_{\pm_2}(t')\psi_0,S_{\pm_3}(t')\psi_0 \Big)  S_{\pm_4}(t')\psi_0 \right](t')dt' 	
\end{eqnarray}
where $S_\pm(t)= e^{-\pm it|D|}$ for $\ell=1,2$. From the $C^3$ smoothness we have that 
\begin{eqnarray}\label{smooth}
\sup_{0\le t\le T}\left\| \sum_{\pm_j,j=1,2,3,4}\int_0^t S_{\pm_1}(t-t')\Pi^{\pm_1}(D) \left[ \mathcal N_\ell \Big(S_{\pm_2}(t')\psi_0, S_{\pm_3}(t')\psi_0 \Big) S_{\pm_4}(t')\psi_0 \right]dt'  \right\|_{H^s} \les \|\psi_0\|_{H^s}^3
\end{eqnarray}
for a local existence time $T$ and $j=1,2$. However we show that \eqref{smooth} fails for $s<s(\ell)$. The explicit statement is as follows:
\begin{prop}\label{prop-illpose}
	Let $\ell=1,2$. Assume that $s< s(\ell)$. Then the estimate
	\begin{eqnarray}\label{illpose}
	\sup_{0\le t\le T}\left\| \mathcal{L}_\ell(\vp)(t)  \right\|_{H^s} \les \|\vp\|_{H^s}^3.
	\end{eqnarray}
	fails to hold for all $\varphi \in H^s$, where $\mathcal{L}_\ell(\vp)(t)=\underset{\pm_j,j=1,\cdots,4}{\sum}\mathcal{L}_\ell^{1,\cdots,4}(\varphi)(t)$ with
	$$
	\mathcal{L}_\ell^{1,\cdots,4}(\varphi)(t) = \int_0^t  S_{\pm_1}(t-t')\Pi^{\pm_1}(D) \mathcal N_\ell\Big( S_{\pm_2}(t')\vp,  S_{\pm_3}(t')\vp,  \beta S_{\pm_4}(t')\vp \Big) dt'.
	$$
\end{prop}

\begin{proof}
	The proof of Proposition \ref{prop-illpose} is proven by contradiction. For this purpose, let us assume that the \eqref{illpose} holds. Fix $ \lam \gg 1$. We first choose $\mu = \lam^{1-\varepsilon}$ for fixed $0<\varepsilon \ll 1$. Let us define a box
	\begin{align*}
	B_{\mu} &= \{ \xi=(\xi_1,\xi_2) : |\xi_1 - \lam|\les \mu,\;|\xi_2|\les \mu \}
	\end{align*}
	and consider $\vp = \left(\begin{array}{c} \mathcal{F}_\xi^{-1}\chi_{B_\mu}\\ 0 \end{array} \right) $. Then we have  $\|\varphi\|_{H^s}  \sim \mu\lam^s$.

	To lead a contradiction we adopt a following estimate:
	\begin{eqnarray}\label{gtr-ine}
	\left|   \sum_{\pm_j,j=1,\cdots,4} \mathcal{F}_x\Big[\mathcal{L}_\ell^{1\cdots4}(\vp)(t)\Big](\xi)     \right| \gtrsim  t\mu^{3+2s(j)}.
	\end{eqnarray}
    We now prove the \eqref{gtr-ine}. By taking Fourier transform we see that
    \begin{align*}
    	& \mathcal{F}_x\Big[\mathcal{L}_\ell^{1\cdots4}(\vp)(t)\Big](\xi)  \\
    	&=  \Pi^{\pm}(\xi)   \int_0^t \int e^{-\pm_1i(t-t')|\xi|} \mathcal{F}_x \left[ \mathcal N_\ell\Big( S_{\pm_2}(t')\vp,  S_{\pm_3}(t')\vp   \Big) \right](\sigma)  \mathcal{F}_x  \left[S_{\pm_4}(t')\vp \right](\xi-\sigma)  d\sigma dt'\\
    	&=-\Pi^{\pm}(\xi)  \int_{|\sigma| \les \mu}\int_{-B_\mu} \mathbf{p}_{1\cdots4}(t,\xi,\sigma,\zeta)|\sigma|^{-1+2s(j)}\chi_{B_\mu}(-\zeta)\chi_{B_\mu}(\sigma - \zeta) \chi_{B_\mu}(\xi-\sigma) d\zeta d\sigma\\
    \end{align*}
    where $-B_\mu := \{ \xi= (\xi_1,\xi_2) : (-\xi_1, -\xi_2) \in B_\mu\}$     and
    \begin{align*}
    	\mathbf{p}_{1\cdots4}(t,\xi,\sigma,\zeta) &:= \int_0^t e^{-i(\pm_1(t-t')|\xi|\pm_2t'|\zeta|\pm_3t'|\sigma-\zeta|\pm_4t'|\sigma|)}dt'\\
    	&= \frac{e^{-\pm_1it|\xi|}(e^{it\omega}-1)}{i\omega}
    \end{align*}
    with
    \begin{align*}
    	\omega = \pm_1 |\xi| \pm_2 |\zeta| \pm_3 |\sigma - \zeta| \pm_4 |\sigma|.
    \end{align*}
    From the support condition it follows that $|\sigma| \les 2\mu$, provided $\xi \in B_{3\mu} $. Then $|\omega| \les \lam$.
 
    We set $t= \delta\lam^{-1-\varepsilon}$ for fixed $0<\delta \ll 1$. Since $|t\omega| \ll 1$ for $\lam$ large enough, we get
    \begin{align*}
    	\sum_{\pm_j,j=1,2,3,4}\mathbf{p}_{1\cdots4}(t,\xi,\sigma,\zeta) &\qquad = \sum_{\pm_j,j=1,2,3,4} te^{-\pm_1it|\xi|}\left( \frac{\cos(t\omega)-1}{it\omega} +i \frac{\sin(t\omega)}{it\omega}\right)\\
    	&\qquad = \sum_{\pm_j,j=1,2,3,4} te^{-\pm_1it|\xi|}(O_\pm(\delta) +i)\\
    	&\qquad = \sum_{\pm_j,j=1,2,3,4} te^{-\pm_1it|\xi|}O_\pm(\delta) +i\sum_{\pm_j,j=1,2,3,4} te^{-\pm_1it|\xi|}\\
    	&\qquad = \sum_{\pm_j,j=1,2,3,4} te^{-\pm_1it|\xi|}O_\pm(\delta) +8it\cos(t|\xi|)\\
    	&\qquad = \sum_{\pm_j,j=1,2,3,4} te^{-\pm_1it|\xi|}O_\pm(\delta) +8it(1+O(\delta)).
    \end{align*}   
    Hence we obtain
    \begin{align*}
    	&\left|   \sum_{\pm_j,j=1,2,3,4}  \mathcal{F}_x\Big[\mathcal{L}_\ell^{1\cdots4}(\vp)(t)\Big](\xi)  \right| \\
    	&\qquad\gtrsim      \left| \sum_{\pm_j,j=1,2,3,4}\int_{|\sigma| \les \mu}\int_{-B_\mu} \mathbf{p}_{1\cdots4}(t,\xi,\sigma,\zeta)|\sigma|^{-1+2s(\ell)}\chi_{B_\mu}(-\zeta)\chi_{B_\mu}(\sigma - \zeta) \chi_{B_\mu}(\xi-\sigma) d\zeta d\sigma\right|\\
    	&\qquad\gtrsim t \left|  \int_{|\sigma| \les \mu}\int_{-B_\mu}|\sigma|^{-1+2s(\ell)}\chi_{B_\mu}(-\zeta)\chi_{B_\mu}(\sigma - \zeta) \chi_{B_\mu}(\xi-\sigma) d\zeta d\sigma\right|\\
    	&\qquad\gtrsim t\mu^{3+2s(\ell)}.
    \end{align*}
    Therefore we get \eqref{gtr-ine}.

    We return to the main proof. Since $\mathcal{F}_x\Big[\mathcal{L}_\ell^{1\cdots4}(\vp)(t)\Big](\xi) =0$ for $\xi \notin B_{3\mu}$, the \eqref{gtr-ine} yields that
	\begin{align*}
	\left\| \mathcal{L}_\ell(\vp)(t)  \right\|_{H^s} = \left\|\langle\xi\rangle^s  \sum_{\pm_j,j=1,2,3,4} \mathcal{F}_x\Big[\mathcal{L}_\ell^{1\cdots4}(\vp)(t)\Big](\xi)  \right\|_{L_{\xi}^2} &\gtrsim t\mu^{3+2s(\ell)}\left\|\langle \xi \rangle^s \right\|_{L_{\xi}^2(B_{3\mu})} \gtrsim t \mu^{4+2s(\ell)}\lam^s.
	\end{align*}
	This gives us that
	\begin{eqnarray}\label{failure2}
	t\mu^{4+2s(\ell)}\lam^{s} \les \left\|\langle\xi\rangle^s  \sum_{\pm_j,j=1,\cdots,4}  \mathcal{F}_x\Big[\mathcal{L}_j^{1\cdots4}(\vp)(t)\Big](\xi) \right\|_{L_{\xi}^2} \les \mu^{3}\lam^{3s}.
	\end{eqnarray}
	
	Therefore, by \eqref{failure2} and $t=\delta\lam^{-1-\varepsilon}$, we have
	\begin{align}\label{contradiction}
	\delta \les \mu^{-1-2s(\ell)} \lam^{2s+1+\ve}=\lam^{2s + 2s(\ell) +2\varepsilon \big(1+s(\ell) \big) }.
	\end{align}
	Then since  the \eqref{contradiction} does not hold for $s<s(\ell)$ and $\lam \gg 1$, we reach a contradiction. This completes the proof of Proposition \ref{prop-illpose}.
\end{proof}


\section*{Acknowledgments}
 This work was supported by National University Promotion Development Project in 2019.



\begin{thebibliography}{00}
	
\bibitem{phymod} J. Arbunich and S. Christof, \textit{Rigorous derivation of nonlinear Dirac equations for wave propagation in honeycomb structures}, Jour. of Math, Phys. \textbf{59.1} (2018), 011509.

\bibitem{bapese} O. Bahat-Treidel, O. Peleg, and M. Segev, \textit{Symmetry breaking in honeycomb photonic lattices}, Opt. Lett. \textbf{33} (2008), 2251--2253.


\bibitem{behe} I. Bejenaru and S. Herr, \textit{The cubic Dirac equation: small initial data in $H^{\frac12}(\mathbb{R}^2)$}, Comm. Math. Phys., \textbf{343} (2016), 515--562.

\bibitem{beta}  I. Bejenaru and T. Tao, \textit{Sharp well-posedness and ill-posedness results for a quadratic non-linear Schr\"odinger
equation}, J. Funct. Anal. \textbf{233} (2006), 228--259.

\bibitem{bourgain} J. Bourgain, \textit{Fourier transform restriction phenomena for certain lattice subsets and applications to nonlinear evolution equations}, I. Schrödinger
equations, Geom. Funct. Anal. \textbf{3(2)} (1993), 107--156.



\bibitem{castro} A. H. Castro Neto, F. Guinea, N. M. R. Peres, K. S. Novoselov, and A. K. Geim, \textit{The electronic properties of graphene}, Rev. Mod. Phys. \textbf{81} (2009), 109--162.

\bibitem{cholee} Y. Cho and K. Lee, \textit{Small data scattering of Dirac equations with Yukawa type potentials in $L_x^2(\mathbb R^2)$}, preprint.

\bibitem{choleeoz} Y. Cho, K. Lee, and T. Ozawa, \textit{Small data scattering of 2d Hatree type Dirac equations}, preprint.

\bibitem{choz}  Y. Cho and T. Ozawa, \textit{On the semirelativistic Hartree-type equation}, SIAM J. Math. Anal. \textbf{38} (2006), 1060-1074.


\bibitem{chozsashim} Y. Cho, T. Ozawa, H. Sasaki, and  Y. Shim  \textit{Remarks on the semirelativistic Hartree equations}, Discrete Contin. Dyn. Syst. \textbf{23} (2009), no. 4, 1277--1294.


\bibitem{din} V. Dinh, \textit{On the Cauchy problem for the nonlinear semi-relativistic equation in Sobolev spaces}, Des. Con. Dyn. Sys. \textbf{38} (2018), 1127--1148

\bibitem{fugeooz} K. Fujiwara, V. Georgiev, and T. Ozawa, \textit{On global well-posedness for nonlinear semirelativistic equations in some scaling subcritical and critical cases} Jour. de math. pure. et appl. \textbf{136} (2020), 239--256.


\bibitem{gtv} J. Ginibre, Y. Tsutsumi, and G. Velo, \textit{On the Cauchy problem for the Zakharov system}, Jour. of Func. Anal. \textbf{151.2} (1997), 384--436.


\bibitem{hame} R. Hajj and F. M\`ehats, \textit{Analysis of models for quantum transport of electrons in graphene layers}, Math. Models Methods Appl. Sci. \textbf{24(11)} (2014), 2287-–2310.



\bibitem{hele} S. Herr and E. Lenzmann, \textit{The Boson star equation with initial data of low regularity}, Nonlinear Analysis \textbf{97} (2014), 125--137.



\bibitem{kapo}  T.  Kato and  G.  Ponce,  \textit{Commutator  estimates  and  the  Euler  and  Navier–Stokes  equations},  Comm.  Pure  Appl. Math. \textbf{41(7)} (1988), 891--907.


\bibitem{kgo} F. Kazumasa, V. Georgiev, and T. Ozawa, \textit{Higher order fractional Leibniz rule}, Journal of Fourier Analysis and Applications. \textbf{24.3} (2018), 650--665.

\bibitem{kpv2}  C.E. Kenig, G. Ponce, and L. Vega, \textit{Well-posedness and scattering results for the generalized Korteweg–de Vries equation via the contraction principle}, Comm. Pure Appl. Math. \textbf{46} (1993), 527--620.


\bibitem{kpv} C.E. Kenig,  G. ponce, and L. Vega, \textit{A bilinear estimate with applications to the KdV equations}, J. Amer. Math. Soc. \textbf{9} (1996), 573--603.  



\bibitem{lee}  K. Lee, \textit{Low regularity  well-posedness of Hartree type Dirac equations in 2,3-dimensions}, preprint.


\bibitem{len} E. Lenzmann, \textit{Well-posedness for semi-relativistic Hartree equations of critical type}, Math. Phys. Anal. Geom. \textbf{10 (1)} (2007), 43–-64.

\bibitem{mst} L. Molinet, Jean-Claude Saut, and N. Tzvetkov, \textit{Ill-Posedness Issues for the Benjamin--Ono and Related Equations}, SIAM journal on mathematical analysis \textbf{33.4} (2001), 982--988.

\bibitem{sel} S. Selberg, \textit{Bilinear Fourier restriction estimates related to the 2D wave equation}, Adv. Diff. Equ. \textbf{16.7/8} (2011), 667--690.

\bibitem{tao}  T. Tao,  \textit{Multilinear weighted convolution of $L^2$ functions, and applications to nonlinear dispersive equations}, American Journal of Mathematics \textbf{123(5)} (2001), 839--908.

\bibitem{tes2d} A. Tesfahun, \textit{Long-time behavior of solutions to cubic Dirac equation with Hartree type nonlinearity in $\mathbb{R}^{1+2}$}, International Mathematics Research Notices \textbf{} (2018), 1--50.


\bibitem{tes3d} \bysame, \textit{Small data scattering for cubic Dirac equation with Hartree type nonlinearity in $\mathbb R^{1+3}$}, SIAM J. Math. Anal. \textbf{52} (2020), no. 3, 2969--3003.

\bibitem{yang} C. Yang, \textit{Scattering results for Dirac Hartree-type equations with small initial data}, Comm. Pure. Appl. Anal. \textbf{18 (4)} (2019), 1711--1734.


\end{thebibliography}
\end{document}